\newcommand{\al}{\alpha}               
\newcommand{\be}{\beta}
\newcommand{\ga}{\gamma}               

\newcommand{\de}{\delta}

\newcommand{\Lb}{\Lambda}
\newcommand{\sig}{\sigma}

\newcommand{\veps}{\varepsilon}        
\newcommand{\vphi}{\varphi}

\newcommand{\cal}{\mathcal}
\newcommand{\cala}{{\cal A}}           
\newcommand{\calb}{{\cal B}}
           
\newcommand{\calf}{{\cal F}}
\newcommand{\calg}{{\cal G}}           

\newcommand{\cali}{{\cal I}}           
           
\newcommand{\calp}{{\cal P}}
\newcommand{\calr}{{\cal R}}           
\newcommand{\cals}{{\cal S}}

\newcommand{\diam}{{\rm diam}}      

\newcommand{\Fix}{{\rm Fix}}

\newcommand{\spec}{{\rm spec}}

\newcommand{\incl}{\subseteq}

\newcommand{\es}{\emptyset}          
\newcommand{\sm}{\setminus}
      
\newcommand{\limpl}{\Longrightarrow}

\newcommand{\oo}{\infty}

\newcommand{\wt}{\widetilde}

\newcommand{\n}{\noindent}

\def\R+oo{R_+\cup\{\oo\}}

\def\dtends   {\stackrel {\it d}{\longrightarrow}}
\def\0dtends  {\stackrel {\it 0d}{\longrightarrow}}

\newcommand{\barr}{\begin{array}}         
\newcommand{\earr}{\end{array}}
\newcommand{\bcor}{\begin{corollary}}    
\newcommand{\ecor}{\end{corollary}}
\newcommand{\ben}{\begin{enumerate}}      
\newcommand{\een}{\end{enumerate}}
\newcommand{\beq}{\begin{equation}}       
\newcommand{\eeq}{\end{equation}}
\newcommand{\bex}{\begin{example}}        
\newcommand{\eex}{\end{example}}
\newcommand{\bit}{\begin{itemize}}        
\newcommand{\eit}{\end{itemize}}
\newcommand{\blemma}{\begin{lemma}}       
\newcommand{\elemma}{\end{lemma}}
\newcommand{\bproof}{\begin{proof}}       
\newcommand{\eproof}{\end{proof}}
\newcommand{\bprop}{\begin{proposition}}  
\newcommand{\eprop}{\end{proposition}}
\newcommand{\brem}{\begin{remark}}        
\newcommand{\erem}{\end{remark}}
\newcommand{\btab}{\begin{tabular}}       
\newcommand{\etab}{\end{tabular}}
\newcommand{\btheorem}{\begin{theorem}}   
\newcommand{\etheorem}{\end{theorem}}

\documentclass[reqno]{amsart}

\newtheorem{theorem}{\bf Theorem}
\newtheorem{corollary}{\bf Corollary}
\newtheorem{example}{\bf Example}
\newtheorem{lemma}{\bf Lemma}
\newtheorem{proposition}{\bf Proposition}
\newtheorem{remark}{\bf Remark}

\begin{document}

\title
[Contractive Maps in LT-Relational Metric Spaces]
{CONTRACTIVE MAPS IN LOCALLY TRANSITIVE \\
RELATIONAL METRIC SPACES}

\author{Mihai Turinici}
\address{
"A. Myller" Mathematical Seminar;
"A. I. Cuza" University;
700506 Ia\c{s}i, Romania
}
\email{mturi@uaic.ro}




\begin{abstract}
Some fixed point results are given for a class of
Meir-Keeler contractive maps acting on 
metric spaces endowed with locally transitive relations. 
Technical  connections with the related statements due to 
Berzig et al
[Abstr. Appl. Anal., Volume 2013, Article ID 259768]
are also being discussed.
\end{abstract}

\maketitle

\section{Introduction}
\setcounter{equation}{0}

Let $X$ be a nonempty set.
Call the subset $Y$ of $X$, 
{\it almost-singleton} (in short: {\it asingleton})
provided $y_1,y_2\in Y$ implies $y_1=y_2$;
and {\it singleton},
if, in addition, $Y$ is nonempty;
note that, in this case,
$Y=\{y\}$, for some $y\in X$. 
Take a {\it metric} 
$d:X\times X\to R_+:=[0,\oo[$ over $X$;
as well as a selfmap 
$T\in \calf(X)$.
[Here, for each couple $A,B$ of nonempty sets,
$\calf(A,B)$ denotes the class of all 
{\it functions} 
from $A$ to $B$; 
when $A=B$, we write $\calf(A)$ in place of $\calf(A,A)$].
Denote $\Fix(T)=\{x\in X; x=Tx\}$;
each point of this set is referred to as 
{\it fixed} under $T$.
Concerning the existence and uniqueness of such points,
a basic result is the 1922 one due to
Banach \cite{banach-1922}.
Call the selfmap $T$, {\it $(d;\al)$-contractive} (where $\al\ge 0$), if
\ben
\item[] (a01)\ \ 
$d(Tx,Ty)\le \al d(x,y)$,\ for all $x,y\in X$.
\een

\btheorem  \label{t1}
Assume that $T$ is $(d;\al)$-contractive, for some $\al\in [0,1[$.
In addition, let $X$ be $d$-complete.
Then,

{\bf i)}
$\Fix(T)$ is a singleton, $\{z\}$;\ \ 
{\bf ii)}
$T^nx \dtends z$ as $n\to \oo$, for each $x\in X$.
\etheorem

This result 
(referred to as: Banach's fixed point theorem) 
found some basic applications to 
the operator equations theory.
As a consequence, a multitude of 
extensions for it were proposed.
Here, we shall be interested in the 
{\it relational} way of
enlarging Theorem \ref{t1},
based on contractive conditions like
\ben
\item[] (a02)\ \ 
$F(d(Tx,Ty),d(x,y),d(x,Tx),d(y,Ty),d(x,Ty),d(y,Tx))\le 0$,\\
for all $x,y\in X$ with $x\calr y$;
\een
where $F:R_+^6\to R$ is a function, 
and $\calr$ is a 
{\it relation} over $X$.
Note that, when 
$\calr$ is the 
{\it trivial} relation
(i.e.: $\calr=X\times X$), 
a large list of such 
contractive maps is provided in 
Rhoades \cite{rhoades-1977}.
Further, when $\calr$ is 
an {\it order} on $X$,
a first result is the 1986 one obtained by
Turinici \cite{turinici-1986},
in the realm of ordered metrizable uniform spaces.
Two decades after, 
this fixed point statement
was re-discovered 
(in the ordered metrical setting)
by
Ran and Reurings \cite{ran-reurings-2004};
see also
Nieto and Rodriguez-Lopez 
\cite{nieto-rodriguez-lopez-2005};
and, since then, the 
number of such results 
increased rapidly.
On the other hand, 
when $\calr$ is an {\it amorphous} relation
over $X$, an appropriate statement of this type 
is the 2012 one due to
Samet and Turinici \cite{samet-turinici-2012}.
The "intermediary" particular case of 
$\calr$ being 
{\it finitely transitive} was 
recently obtained by
Karapinar and Berzig \cite{karapinar-berzig-2013},
under a class of $(\al\psi,\be\vphi)$-contractive 
conditions suggested by
Popescu \cite{popescu-2011}.
It is our aim in the following to 
give further extensions of these results, when

{\bf i)} 
the contractive conditions 
are taken after the model in 
Meir and Keeler \cite{meir-keeler-1969}

{\bf ii)}
the finite transitivity of $\calr$ is being
assured in a "local" way. 

\n
Further aspects will be delineated elsewhere.

\section{Preliminaries}
\setcounter{equation}{0}

Throughout this exposition, 
the ambient axiomatic system is 
Zermelo-Fraenkel's (abbreviated: (ZF)). 
In fact, the {\it reduced} system
(ZF-AC) will suffice;
here, (AC) stands for the {\it Axiom of Choice}.
The notations and basic facts 
to be used in this reduced system
are standard.
Some important ones are described below.

{\bf (A)}
Let $X$ be a nonempty set.
By a {\it relation} over $X$, we mean
any nonempty part $\calr\incl X\times X$.
For simplicity, we sometimes write
$(x,y)\in \calr$ as $x\calr y$.
Note that $\calr$ may be regarded 
as a mapping between $X$ and
$\calp(X)$ (=the class of all subsets in $X$).
In fact, denote for $x\in X$:
$X(x,\calr)=\{y\in X; x\calr y\}$
(the {\it section} of $\calr$ through $x$);
then, the desired mapping representation is 
[$\calr(x)=X(x,\calr)$, $x\in X$].

Among the classes of relations to be used,
the following ones 
(listed in an "increasing" scale)
are important for us:
\ben
\item[] (P0)\ \ 
$\calr$ is {\it amorphous}; i.e.:
it has no specific properties at all
\item[] (P1)\ \ 
$\calr$ is an {\it order};
i.e.: it is
{\it reflexive}
[$x\calr x$, $\forall x\in X$],
{\it transitive}
[$x\calr y$ and $y\calr z$ imply $x\calr z$]
and
{\it antisymmetric}
[$x\calr y$ and $y\calr x$ imply $x=y$]
\item[] (P2)\ \ 
$\calr$ is a {\it quasi-order};
i.e.: it is reflexive and transitive
\item[] (P3)\ \ 
$\calr$ is transitive (see above).
\een

A basic ordered structure is $(N,\le)$;
here, $N=\{0,1,...\}$ is the set of natural numbers and
$(\le)$ is defined as
$m\le n$ iff $m+p=n$, for some $p\in N$.
For each $n\in N(1,\le)$, let $N(n,>):=\{0,...,n-1\}$ stand for
the {\it initial interval} (in $N$) induced by $n$.
Any set $P$ with 
$P\sim N$ (in the sense: there exists a bijection from 
$P$ to $N$) will be referred to as 
{\it effectively denumerable}. 
In addition, given some 
natural number $n\ge 1$, any set $Q$ with 
$Q\sim N(n,>)$ will be said to be {\it $n$-finite};
when $n$ is generic here, we say that $Q$ is {\it finite}.
Finally, the (nonempty) set $Y$ is called
(at most) {\it denumerable} iff it is either 
effectively denumerable or finite.

Given the relations $\calr$, $\cals$ over $X$,
define their {\it product}
$\calr\circ \cals$ as
\ben
\item[] (b01)\ \ 
$(x,z)\in \calr\circ \cals$
if, there exists $y\in X$ with
$(x,y)\in \calr$, $(y,z)\in \cals$.
\een
This allows us 
to introduce the powers of a relation $\calr$ as
\ben
\item[] (b02)\ \ 
$\calr^0=\cali$,\ 
$\calr^{n+1}=\calr^n \circ \calr$,\ $n\in N$.
\een
[Here, $\cali=\{(x,x); x\in X\}$ is the 
{\it identical relation} over $X$].
The following properties will be useful in the sequel:
\beq \label{201}
\calr^{m+n}=\calr^m \circ \calr^n,\
(\calr^m)^n=\calr^{mn},\ \forall m,n\in N.
\eeq

Given $k\in N(2,\le)$, let us say that 
$\calr$ is {\it $k$-transitive},
if $\calr^k\incl \calr$;
clearly, {\it transitive} is identical with 
{\it 2-transitive}.
We may now complete the increasing scale above as 
\ben
\item[] (P4)\ \ 
$\calr$ is {\it finitely transitive};
i.e.:
$\calr$ is $k$-transitive 
for some $k\ge 2$
\item[] (P5)\ \ 
$\calr$ is {\it locally finitely transitive};
i.e.:
for each 
(effectively) denumerable subset $Y$ of $X$,
there exists $k=k(Y)\ge 2$, such that 
the restriction to $Y$ of $\calr$ is 
$k$-transitive
\item[] (P6)\ \ 
$\calr$ is {\it trivial};
i.e.: $\calr=X\times X$;
hence,
[$x\calr y$, $\forall x,y\in X$].
\een  

Concerning these concepts,
the following property will be useful.
Call the sequence $(z_n; n\ge 0)$ in $X$, 
{\it $\calr$-ascending}, if $z_i\calr z_{i+1}$ for all $i\ge 0$.

\blemma \label{le1}
Let the $\calr$-ascending sequence 
$(z_n; n\ge 0)$ in $X$ 
and the natural number $k\ge 2$, be such that
\ben
\item[] (b03)\ \ 
$\calr$ is $k$-transitive on $Z:=\{z_n; n\ge 0\}$.
\een
Then, necessarily, 
\beq \label{202}
(\forall r\ge 0):\ 
[(z_i,z_{i+1+r(k-1)})\in \calr,\ \forall i\ge 0].
\eeq
\elemma

\bproof
We shall use the induction with respect to $r$.
First, by the choice of our sequence,
$(z_i,z_{i+1})\in \calr$;
whence, the case $r=0$ holds.
Moreover, by definition,
$(z_i,z_{i+k})\in \calr^k$;
and this, along with 
the $k$-transitive property, gives
$(z_i,z_{i+k})\in \calr$;
hence, the case of $r=1$ holds too.
Suppose that this property holds for some
$r\ge 1$;
we claim that it holds as well for $r+1$.
In fact, let $i\ge 0$ be arbitrary fixed.
Again by the choice of our sequence,
$(z_{i+1+r(k-1)},z_{i+1+(r+1)(k-1)})\in \calr^{k-1}$;
so that, by the 
inductive hypothesis (and properties of relational product)
$$
(z_i,z_{i+1+(r+1)(k-1)})\in \calr \circ \calr^{k-1}=\calr^k;
$$
and this, along with the $k$-transitive condition, gives
$(z_i,z_{i+1+(r+1)(k-1)})\in \calr$.
The proof is thereby complete.
\eproof

{\bf (B)}
Let $(X,d)$ be a metric space.
We introduce a $d$-convergence and $d$-Cauchy structure on $X$ 
as follows.
By a {\it sequence} in $X$, we mean any mapping $x:N\to X$.
For simplicity reasons, 
it will be useful to denote it as 
$(x(n); n\ge 0)$, or $(x_n; n\ge 0)$;
moreover, when no confusion can arise, 
we further simplify this notation as 
$(x(n))$ or $(x_n)$, respectively.
Also, any sequence $(y_n:=x_{i(n)}; n\ge 0)$ with
$i(n)\to \oo$ as $n\to \oo$ will be referred to as a
{\it subsequence} of $(x_n; n\ge 0)$.
Given the sequence $(x_n)$ in $X$ and the point $x\in X$,
we say that $(x_n)$, {\it $d$-converges} to $x$ 
(written as: $x_n \dtends x$)
provided $d(x_n,x)\to 0$ as $n\to \oo$; i.e.,
\ben
\item[]  
$\forall\veps> 0$, $\exists i=i(\veps)$:\ \ 
$i\le n \limpl d(x_n,x)< \veps$.
\een
The set of all such points $x$ will be denoted 
$\lim_n (x_n)$; note that, it is an asingleton, 
because $d$ is triangular symmetric;
if $\lim_n (x_n)$ is nonempty, then 
$(x_n)$ is called {\it $d$-convergent}.
We stress that the introduced convergence concept 
$(\dtends)$ does match the standard
requirements in
Kasahara \cite{kasahara-1976}.
Further, call the sequence $(x_n)$, {\it $d$-Cauchy} 
when $d(x_m,x_n)\to 0$ as $m,n\to \oo$, $m< n$; i.e.,
\ben
\item[]  
$\forall\veps> 0$, $\exists j=j(\veps)$:\ \ 
$j\le m< n \limpl d(x_m,x_n)< \veps$.
\een
As $d$ is triangular symmetric, 
any $d$-convergent sequence is $d$-Cauchy too;
but, the reciprocal is not in general true.
Concerning this aspect, 
note that any $d$-Cauchy sequence $(x_n; n\ge 0)$ 
is {\it $d$-semi-Cauchy}; i.e.,
\ben
\item[]  
$d(x_n,x_{n+1})\to 0$\ 
(hence, $d(x_n,x_{n+i})\to 0$, $\forall i\ge 1$), as $n\to \oo$.
\een
But, the reciprocal is not in general true.

The introduced concepts allow us
to give a useful property:

\blemma \label{le2}
The mapping $(x,y)\mapsto d(x,y)$ is $d$-Lipschitz, 
in the sense  
\beq \label{203}
|d(x,y)-d(u,v)|\le d(x,u)+d(y,v),\ 
\forall (x,y),(u,v)\in X\times X.
\eeq
As a consequence, this map is $d$-continuous; i.e.,
\beq \label{204}
\mbox{
$x_n\dtends x$, $y_n\dtends y$ 
imply $d(x_n,y_n)\to d(x,y)$.
}
\eeq
\elemma

The proof is immediate, by the usual properties of 
the ambient metric $d(.,.)$;
we do not give details.

{\bf (C)}
Let $(X,d)$ be a metric space;
and $\calr\incl X\times X$ be a
(nonempty) relation over $X$;
the triple 
$(X,d,\calr)$ will be referred to as a
{\it relational metric space}.
Further, take some $T\in \calf(X)$.
Call the subset $Y$ of $X$, {\it $\calr$-almost-singleton}
(in short: {\it $\calr$-asingleton}) provided
$y_1,y_2\in Y$, $y_1\calr y_2$ $\limpl$ $y_1=y_2$;
and {\it $\calr$-singleton} when, 
in addition, $Y$ is nonempty.
We have to determine 
circumstances under which $\Fix(T)$ be nonempty;
and, if this holds, to establish 
whether $T$ is {\it fix-$\calr$-asingleton} 
(i.e.: $\Fix(T)$ is $\calr$-asingleton);
or, equivalently:
$T$ is {\it fix-$\calr$-singleton} 
(in the sense: $\Fix(T)$ is $\calr$-singleton);
To do this, we start from the working hypotheses
\ben
\item[] (b04)\ \ 
$T$ is $\calr$-semi-progressive: 
$X(T,\calr):=\{x\in X; x\calr Tx\}\ne \es$ 
\item[] (b05)\ \ 
$T$ is $\calr$-increasing:\ $x\calr y$ implies $Tx\calr Ty$.
\een

The basic 
directions under which the investigations be conducted 
are described by the list below,
comparable with the one in
Turinici \cite{turinici-2011}:

{\bf 2a)}
We say that $T$ is a 
{\it Picard operator} (modulo $(d,\calr)$)
if, for each $x\in X(T,\calr)$,
$(T^nx; n\ge 0)$ is $d$-convergent

{\bf 2b)}
We say that $T$ is a
{\it strong Picard operator} (modulo $(d,\calr)$)
when, for each $x\in X(T,\calr)$,
$(T^nx; n\ge 0)$ is $d$-convergent
and $\lim_n(T^nx)\in \Fix(T)$

{\bf 2c)}
We say that $T$ is a 
{\it globally strong Picard operator} 
(modulo $(d,\calr)$)
when it is a strong Picard operator 
(modulo $(d,\calr)$) and 
$T$ is fix-$\calr$-asingleton (hence, fix-$\calr$-singleton).

The sufficient (regularity) conditions 
for such properties are being
founded on {\it ascending orbital} concepts
(in short: (a-o)-concepts).
Remember that the sequence $(z_n; n\ge 0)$ in $X$ is called  
{\it $\calr$-ascending}, if $z_i\calr z_{i+1}$ for all $i\ge 0$;
further, let us say that $(z_n; n\ge 0)$ is 
{\it $T$-orbital}, when 
it is a subsequence of $(T^n x; n\ge 0)$, 
for some $x\in X$;
the intersection of these notions is just the precise one.

{\bf 2d)}
Call $X$, {\it $(a-o,d)$-complete}, provided
(for each (a-o)-sequence) $d$-Cauchy 
$\limpl$ $d$-convergent

{\bf 2e)}
We say that $T$ is {\it $(a-o,d)$-continuous}, if
($(z_n)$=(a-o)-sequence and $z_n\dtends z$) imply  
$Tz_n\dtends Tz$

{\bf 2f)}
Call $\calr$, {\it $(a-o,d)$-almost-selfclosed}, if:
whenever the (a-o)-sequence $(z_n; n\ge 0)$ in $X$
and the point $z\in X$ fulfill $z_n\dtends z$,
there exists a subsequence $(w_n:=z_{i(n)}; n\ge 0)$ of
$(z_n; n\ge 0)$ with 
$w_n\calr z$, for all $n\ge 0$.

When the orbital properties are ignored,
these conventions give us 
{\it ascending} notions (in short: a-notions).
On the other hand, when the ascending properties
are ignored, the same conventions give us
{\it orbital} notions (in short: o-notions).
The list of these is obtainable
from the previous one; so, further details are
not needed. 
Finally, when $\calr=X\times X$,
the list of such notions is comparable with the one in
Rus \cite[Ch 2, Sect 2.2]{rus-2001}:
because, in this case, $X(T,\calr)=X$.

\section{Meir-Keeler contractions}
\setcounter{equation}{0}

Let $(X,d,\calr)$ be a relational metric space;
and $T$ be a selfmap of $X$;
supposed to be 
$\calr$-semi-progressive
and $\calr$-increasing.
The basic directions and sufficient regularity 
conditions under which 
the problem of determining the fixed points
of $T$ to be solved were already listed.
As a completion of them, we must 
formulate the specific 
metrical contractive conditions upon our data.
These, essentially, consist in a 
"relational" variant of the
Meir-Keeler condition \cite{meir-keeler-1969}.
Assume that
\ben
\item[] (c01)\ \ 
$\calr$ is {\it non-identical}:\ 
[$\wt \calr:=\calr\sm \cali$ is nonempty].
\een
Note that, by definition, the introduced relation writes
\ben
\item[] (c02)\ \  
$x\wt \calr y$\ iff\ [$x\calr y$ and $x\ne y$];
\een
so,  
$\wt \calr$ is {\it irreflexive}
[$x\wt \calr x$ is false, for each $x\in X$].
Denote, for $x,y\in X$:
\ben
\item[] (c03)\ \  
$A_1(x,y)=d(x,y)$,\ $B_1(x,y)=\diam\{x,Tx,y,Ty\}$, \\
$A_2(x,y)=(1/2)[d(x,Tx)+d(y,Ty)]$,\\
$A_3(x,y)=\max\{d(x,Tx),d(y,Ty)\}$, \\
$A_4(x,y)=(1/2)[d(x,Ty)+d(Tx,y)]$.
\een
Then, let us introduce the functions
\ben
\item[] (c04)\ \ 
$B_2=\max\{A_1,A_2\}$, $B_3=\max\{A_1,A_3\}$, 
$B_4=\max\{A_1,A_4\}$,\\
$C_1=\max\{A_1,A_2,A_4\}$, $C_2=\max\{A_1,A_3,A_4\}$, \\
$\calg=\{A_1,B_2,B_3,B_4,C_1,C_2\}$, 
$\calg_1=\{A_1,B_2,B_4,C_1\}$, $\calg_2=\{B_3,C_2\}$.
\een
Note that, for each $G\in \calg$, we have
\beq \label{301}
A_1(x,y)\le G(x,y)\le B_1(x,y),\ \forall x,y\in X.
\eeq
The former of these will be referred to as:
$G$ is {\it sufficient};
note that, by the properties of $d$, we must have
\beq \label{302}
x,y\in X, x\wt \calr y \limpl G(x,y)> 0.
\eeq
And, the latter of these means that:
$G$ is {\it diameter bounded}.

Given $G\in \calg$, we say that $T$ is 
{\it Meir-Keeler $(d,\calr;G)$-contractive}, if
\ben 
\item[] (c05)\ \ 
$x\wt \calr y$ implies $d(Tx,Ty)< G(x,y)$; \\
expressed as:\ $T$ is strictly $(d,\calr;G)$-nonexpansive
\item[] (c06)\ \ 
$\forall \veps> 0$, $\exists \de> 0$:\ 
[$x\wt \calr y$, $\veps< G(x,y)< \veps+\de$] $\limpl$ 
$d(Tx,Ty)\le \veps$; \\
expressed as:\ $T$ has the Meir-Keeler property
(modulo $(d,\calr;G)$).
\een
Note that, by the 
former of these, 
the Meir-Keeler property may be written as
\ben 
\item[] (c07)\ \ 
$\forall \veps> 0$, $\exists \de> 0$:\ 
[$x\wt \calr y$, $G(x,y)< \veps+\de$] $\limpl$ 
$d(Tx,Ty)\le \veps$.
\een

In the following, two basic examples of 
such contractions will be given.

{\bf (A)}
Let $\calf(re)(R_+)$ stand for the class of all
$\vphi\in \calf(R_+)$ with the (strong) {\it regressive} property:
[$\vphi(0)=0$; $\vphi(t)< t$, $\forall t> 0$].
We say that $\vphi\in \calf(re)(R_+)$ 
is {\it Meir-Keeler admissible}, if 
\ben 
\item[] (c08)\ \ 
$\forall \ga> 0$,\ $\exists \be\in ]0,\ga[$,\ $(\forall t)$:\  
$\ga\le t< \ga+\be \limpl \vphi(t)\le \ga$; \\
or, equivalently:\
$\forall \ga> 0$,\ $\exists \be\in ]0,\ga[$,\ $(\forall t)$:\  
$0\le t< \ga+\be \limpl \vphi(t)\le \ga$.
\een 
Now, given 
$G\in \calg$,
$\vphi\in \calf(R_+)$, 
call $T$, 
{\it $(d,\calr;G,\vphi)$-contractive}, if
\ben
\item[] (c09)\ \ 
$d(Tx,Ty)\le \vphi(G(x,y))$, $\forall x,y\in X$, $x\wt \calr y$.
\een

\blemma \label{le3}
Assume that $T$ is 
$(d,\calr;G,\vphi)$-contractive, 
where $\vphi\in \calf(re)(R_+)$ is
Meir-Keeler admissible.
Then, $T$ is Meir-Keeler $(d,\calr;G)$-contractive.
\elemma

\bproof
i)\
Let $x,y\in X$ be such that $x\wt \calr y$.
The contractive condition, and regressiveness of $\vphi$, yield 
$d(Tx,Ty)< G(x,y)$;
so that, 
$T$ is strictly $(d,\calr;G)$-nonexpansive.

ii)\
Let $\veps> 0$ be arbitrary fixed;
and $\de\in ]0,\veps[$ be the number 
assured by the Meir-Keeler admissible property of $\vphi$.
Further, let $x,y\in X$ be such that $x\wt \calr y$ and
$\veps< G(x,y)< \veps+\de$.
By the 
contractive condition and
admissible property,
$$
d(Tx,Ty)\le \vphi(G(x,y))\le \veps;
$$
so that, 
$T$ has the Meir-Keeler property
(modulo $(d,\calr;G)$).

\eproof

Some important classes of such functions
are given below.

{\bf I)}
For any $\vphi\in \calf(re)(R_+)$ and any $s\in R_+^0:=]0,\oo[$, put
\ben
\item[] (c10)\ \ 
$\Lb_+\vphi(s)=\inf_{\veps> 0} \Phi(s+)(\veps)$;\ where
$\Phi(s+)(\veps)=\sup \vphi(]s,s+\veps[)$; 
\item[] (c11)\ \ 
$\Lb^+\vphi(s)=\max\{\vphi(s), \Lb_+\vphi(s)\}$.
\een
By this very definition, we have the representation 
(for all $s\in R_+^0$)
\beq \label{303}
\barr{l}
\Lb^+\vphi(s)=\inf_{\veps> 0} \Phi[s+](\veps);\ 
\mbox{where}\ 
\Phi[s+](\veps)=\sup\{\vphi([s,s+\veps[). 
\earr
\eeq	
From the regressive property of $\vphi$, these limit 
quantities are finite; precisely,
\beq \label{304}
0\le \vphi(s)\le \Lb^+\vphi(s)\le s,\ \ \forall s\in R_+^0.
\eeq

Call $\vphi\in \calf(re)(R_+)$,
{\it Boyd-Wong admissible}, if 
\ben
\item[] (c12)\ \ 
$\Lb^+\vphi(s)< s$\ (or, equivalently: $\Lb_+\vphi(s)< s$),
\ for all $s> 0$.
\een
(This convention is related to the developments in
Boyd and Wong \cite{boyd-wong-1969};
we do not give details).
In particular, $\vphi\in \calf(re)(R_+)$ 
is Boyd-Wong admissible provided 
it is upper semicontinuous at the right on $R_+^0$: 
\ben
\item[]  
$\Lb^+ \vphi(s)= \vphi(s)$,\ 
(or, equivalently:\ $\Lb_+ \vphi(s)\le \vphi(s)$),\ $\forall s\in R_+^0$.
\een
Note that this is fulfilled when 
$\vphi$ is continuous at the right on $R_+^0$;
for, in such a case, 
$\Lb_+ \vphi(s)=\vphi(s)$, $\forall s\in R_+^0$.

{\bf II)}
Call $\vphi\in \calf(re)(R_+)$,
{\it Matkowski admissible}
\cite{matkowski-1975}, 
provided
\ben
\item[] (c13)\ \ 
$\vphi$ is increasing and $\vphi^n(t)\to 0$ as $n\to \infty$, 
for all $t> 0$.
\een
[Here, $\vphi^n$ stands for the $n$-th iterate of $\vphi$].
Note that the obtained class of functions is 
distinct from the above introduced one,
as simple examples show. 

Now, let us say that $\vphi\in \calf(re)(R_+)$ is
{\it Boyd-Wong-Matkowski admissible} 
(abbreviated: BWM-admissible) 
if it is either Boyd-Wong admissible or Matkowski admissible.
The following auxiliary fact will be useful:

\blemma \label{le4}
Let $\vphi\in \calf(re)(R_+)$ be a BWM-admissible function.
Then, $\vphi$ is Meir-Keeler admissible (see above).
\elemma

\bproof 
{\bf (Sketch)}
The former of these is an immediate consequence of
definition.
And, the second one is to be found in
Jachymski \cite{jachymski-1994}.
\eproof

{\bf (B)}
Let us say that  
$(\psi,\vphi)$ is a 
{\it pair of weak generalized altering functions}
in $\calf(R_+)$, if
\ben 
\item[] (c14)\ \ 
$\psi$ is increasing, and 
[$\vphi(0)=0$;\
$\vphi(\veps)> \psi(\veps)-\psi(\veps-0)$, $\forall \veps> 0$] 
\item[] (c15)\ \ 
($\forall \veps> 0$):\ 
$\limsup_n \vphi(t_n)> \psi(\veps+0)-\psi(\veps)$,
whenever $t_n\to \veps++$.
\een

Here, given the sequence $(r_n; n\ge 0)$ in $R$ and the point
$r\in R$, we denoted
\ben
\item[]
$r_n\to r+$ (respectively, $r_n\to r++$), if $r_n\to r$ and \\
$r_n\ge r$ (respectively, $r_n> r$), for all $n\ge 0$ large enough.
\een

Given $G\in \calg$ 
and the couple $(\psi,\vphi)$
of functions in $\calf(R_+)$,
let us say that $T$ is
{\it $(d,\calr;G,(\psi,\vphi))$-contractive},
provided
\ben
\item[] (c16)\ \ 
$\psi(d(Tx,Ty))\le \psi(G(x,y))-\vphi(G(x,y))$,\
$\forall x,y\in X$, $x\wt \calr y$.
\een

\blemma \label{le5}
Suppose that  
$T$ is
$(d,\calr;G,(\psi,\vphi))$-contractive,
for a pair $(\psi,\vphi)$ of
weak generalized altering functions in $\calf(R_+)$.
Then, $T$ is Meir-Keeler 
$(d,\calr;G)$-contractive (see above).
\elemma

\bproof
i)\
Let $x,y\in X$ be such that $x\wt \calr y$.
Then (as $G$ is sufficient), $G(x,y)> 0$; 
so that (by the choice of our pair), 
$\vphi(G(x,y))> 0$; wherefrom
$\psi(d(Tx,Ty))< \psi(G(x,y))$.
This, via [$\psi$=increasing],  
yields $d(Tx,Ty)< G(x,y)$;
so that, 
$T$ is strictly $(d,\calr;G)$-nonexpansive.

ii)\
Assume by contradiction that 
$T$ does not have the Meir-Keeler property
(modulo $(d,\calr;G)$);
i.e., for some $\veps> 0$, 
$$
\forall \de> 0, \exists (x_\de, y_\de)\in \wt \calr:\ 
[\veps< G(x_\de,y_\de)< \veps+\de,\ d(Tx_\de,Ty_\de)> \veps].
$$
Taking a zero converging sequence $(\de_n)$ 
in $R_+^0$,
we get a couple of sequences 
$(x_n; n\ge 0)$ and $(y_n; n\ge 0)$ in $X$, 
so as 
\beq \label{305}
(\forall n):\ \ 
x_n\wt \calr y_n,\ \veps< G(x_n,y_n)< \veps+\de_n,\ d(Tx_n,Ty_n)> \veps.
\eeq
By the contractive condition (and $\psi$=increasing), we get
$$
\psi(\veps)\le \psi(G(x_n,y_n))-\vphi(G(x_n,y_n)),\ \  \forall n;
$$ 
or, equivalently,
\beq \label{306}
(0<)\ \vphi(G(x_n,y_n))\le \psi(G(x_n,y_n))-\psi(\veps),\ \  \forall n.
\eeq 
By (\ref{305}), $G(x_n,y_n)\to \veps++$;
so that, passing to $\limsup$ as $n\to \oo$,
$$
\limsup_n \vphi(G(x_n,y_n))\le \psi(\veps+0)-\psi(\veps).
$$
But, from the hypothesis about $(\psi,\vphi)$,
these relations are contradictory.
This ends the argument.
\eproof

\section{Main result}
\setcounter{equation}{0}

Let $(X,d,\calr)$ be a 
relational metric space.
Further, let 
$T$ be a selfmap of $X$;
supposed to be 
$\calr$-semi-progressive
and $\calr$-increasing.
The basic directions and regularity 
conditions under which 
the problem of determining the fixed points
of $T$ is to be solved, were already listed;
and the contractive type framework was settled.
It remains now to precise the
regularity conditions upon $\calr$.
Denote, for each $x\in X(T,\calr)$,
\ben
\item[] 
$\spec(x)=\{i\in N(1,\le); x\calr T^ix\}$\
(the {\it spectrum} of $x$).
\een
Clearly, $1\in \spec(x)$;
but, the possibility of $\spec(x)=\{1\}$ 
cannot be removed.
This fact remains valid even if
$x\in X(T,\calr)$ is 
{\it orbital admissible}, in the sense:
[$i\ne j$ implies $T^ix\ne T^jx$];
when
the associated orbit 
$T^Nx:=\{T^nx; n\ge 0\}$ 
is effectively denumerable.
But, for the developments below, 
it is necessary that these spectral 
subsets of $N$ should have a finite 
Hausdorff-Pompeiu distance to $N$;
hence, in particular, these must be infinite.  
Precisely, given $k\ge 1$, let us say that 
$\calr$ is {\it $k$-semi-recurrent}
at the orbital admissible $x\in X(T,\calr)$, if:
\ben
\item[]
for each $n\in N(1,\le)$, there exists $q\in \spec(x)$
such that $q\le n< q+k$.
\een
A global version of this convention is the following:
call $\calr$, 
{\it finitely semi-recurrent}
if,
for each orbital admissible $x\in X(T,\le)$,
there exists $k(x)\ge 1$, such that
$\calr$ is $k(x)$-semi-recurrent at $x$.

Assume in the following that
\ben 
\item[] (d01)\ \ 
$\calr$ is finitely semi-recurrent and non-identical.
\een

Our main result in this exposition is

\btheorem \label{t2}
Assume that $T$ is 
Meir-Keeler $(d,\calr;G)$-contractive,
for some $G\in \calg$.
In addition, 
let 
$X$ be $(a-o,d)$-complete;
and one of the following conditions holds

{\bf i)}
$T$ is $(a-o,d)$-continuous

{\bf ii)}
$\calr$ is $(a-o,d)$-almost-selfclosed
and $G\in \calg_1$

{\bf iii)}
$\calr$ is $(a-o,d)$-almost-selfclosed
and $T$ is 
$(d,\calr;G,\vphi)$-contractive,
for a certain Meir-Keeler
admissible function $\vphi\in \calf(re)(R_+)$

{\bf iv)}
$\calr$ is $(a-o,d)$-almost-selfclosed
and $T$ is 
$(d,\calr;G,(\psi,\vphi))$-contractive,
for a certain pair $(\psi,\vphi)$ of
weak generalized altering functions in $\calf(R_+)$.

\n
Then $T$ is a globally strong Picard operator 
(modulo $(d,\calr)$).
\etheorem

\bproof
First, we check the fix-$\calr$-asingleton property.
Let $z_1,z_2\in \Fix(T)$ be such that 
$z_1\calr z_2$; and 
assume by contradiction that $z_1\ne z_2$;
whence $z_1 \wt \calr z_2$.
From the very definitions above,
$$
A_1(z_1,z_2)=A_4(z_2,z_2)=d(z_1,z_2),\
A_2(z_1,z_2)=A_3(z_2,z_2)=0;
$$
whence: $G(z_1,z_2)=d(z_1,z_2)$.
This, via $T$ being  
strictly $(d,\calr;G)$-nonexpansive,
yields an evaluation like
$$
d(z_1,z_2)=d(Tz_1,Tz_2)< G(z_1,z_2);
$$
contradiction; hence the claim.
It remains now to establish the 
strong Picard property (modulo $(d,\calr)$).
The argument will be divided into several steps.

{\bf Part 1.}
We firstly assert that
\beq \label{401}
G(x,Tx)=d(x,Tx),\ 
\mbox{whenever}\  x\wt \calr Tx.
\eeq
Let $x\in X$ be such that 
$x\wt \calr Tx$.
As $T$ is strictly $(d,\calr;G)$-nonexpansive, 
one has
$d(Tx,T^2x)< G(x,Tx)$.
On the other hand, 
$$
\barr{l}
A_4(x,Tx)=(1/2)d(x,T^2x)\le 
(1/2)[d(x,Tx)+d(Tx,T^2x)] \\
=A_2(x,Tx)\le \max\{d(x,Tx),d(Tx,T^2x)\}=A_3(x,Tx).
\earr
$$
This, along with
$$
\barr{l}
d(Tx,T^2x)< A_3(x,Tx) \limpl d(Tx,T^2x)< d(x,Tx) \\
\limpl A_3(x,Tx)=d(x,Tx),
\earr
$$
gives the desired fact.

{\bf Part 2.}
Take some $x_0\in X$;
and put $(x_n=T^nx_0; n\ge 0)$.
If $x_n=x_{n+1}$ for some $n\ge 0$,
we are done;
so, without loss, one may assume that, for each $n\ge 0$,
\ben
\item[] (d02)\ \ 
$x_n\ne x_{n+1}$;\ 
hence, $x_n\wt \calr x_{n+1}$,\ $\rho_n:=d(x_n,x_{n+1})> 0$.
\een
From the preceding part, we derive 
$$ 
\rho_{n+1}=d(Tx_n,Tx_{n+1})< G(x_n,x_{n+1})=\rho_n,\ \forall n;
$$
so that, the sequence $(\rho_n; n\ge 0)$ is 
strictly descending.
As a consequence, $\rho:=\lim_n \rho_n$ exists
as an element of $R_+$.
Assume by contradiction that $\rho> 0$;
and let $\de> 0$ be the number
given by the Meir-Keeler 
$(d,\calr;G)$-contractive condition upon $T$.
By definition, there exists a rank $n(\de)$ such that 
$n\ge n(\de)$ implies $\rho< \rho_n< \rho+\de$;
hence (by a previous representation)
$\rho< G(x_n,x_{n+1})=\rho_n< \rho+\de$.
This, by the Meir-Keeler 
contractive condition we just quoted, yields 
(for the same $n$),
$\rho_{n+1}=d(Tx_n,Tx_{n+1})\le \rho$;
contradiction. Hence, $\rho=0$; so that, 
\beq \label{402}
\rho_n:=d(x_n,x_{n+1})=d(x_n,Tx_n)\to 0,\ 
\mbox{as}\ n\to \oo;
\eeq
i.e. (see above):\ $(x_n; n\ge 0)$ is $d$-semi-Cauchy.

{\bf Part 3.}
Suppose that
\ben
\item[] (d03)\ \ 
there exist $i,j\in N$ such that $i< j$, $x_i=x_j$.
\een
Denoting $p=j-i$, we thus have $p> 0$ and $x_i=x_{i+p}$;
so that
$$
\mbox{
$x_i=x_{i+np}$, $x_{i+1}=x_{i+np+1}$,\ \  for all $n\ge 0$.  
}
$$
By the introduced notations, 
$\rho_i=\rho_{i+np}$, for all  $n\ge 0$.
This, along with 
$\rho_{i+np}\to 0$ as $n\to \infty$,
yields
$\rho_i=0$; 
in contradiction with the initial choice of $(\rho_n; n\ge 0)$.
Hence, our working hypothesis cannot hold; wherefrom
\beq \label{403}
\mbox{
for all $i,j\in N$:\ \ $i\ne j$ implies $x_i\ne x_j$.
}
\eeq

{\bf Part 4.}
As a consequence of this, 
the map $i\mapsto x_i:=T^ix_0$ is injective;
hence, $x_0$ is orbital admissible.
Let $k:=k(x_0)\ge 1$ be the semi-recurrence constant 
of $\calr$ at $x_0$
(assured by the choice of this relation). 
Further, let $\veps> 0$ be arbitrary fixed; 
and $\de> 0$ be the number associated by the
Meir-Keeler $(d,\calr;G)$-contractive property;
without loss, one may assume that $\de< \veps$.
By the $d$-semi-Cauchy property 
and triangular inequality, there exists a rank
$n(\de)\ge 0$, such that 
\beq \label{404}
\barr{l}
(\forall n\ge n(\de)):\ d(x_n,x_{n+1})< \de/4k;\
\mbox{whence} \\
d(x_n, x_{n+h})< h\de/4k \le \de/2,\ \forall h\in \{1,...,2k\}.
\earr
\eeq
We claim that the following relation holds
\beq \label{405}
(\forall s\ge 1):\  
[d(x_n,x_{n+s})< \veps+\de/2,\ \forall n\ge n(\de)];
\eeq
wherefrom, $(x_n; n\ge 0)$ is $d$-Cauchy.
To do this, an induction argument upon $s\ge 1$ 
will be used.
The case $s\in \{1,...,2k\}$ is evident, by the preceding 
evaluation. 
Assume that it holds for 
all $s\in \{1,...,p\}$, where $p\ge 2k$;
we must establish its validity for $s=p+1$.
As $\calr$ is $k$-semi-recurrent at $x_0$,
there exists $q\in \spec(x_0)$ such that 
$q\le p< q+k$;
note that, the former of these yields 
(from the $\calr$-increasing property of $T$),
$x_n\wt \calr x_{n+q}$.
Now, by the inductive hypothesis
and (\ref{404}), 
$$
\barr{l}
d(x_n,x_{n+q}), d(x_{n+1},x_{n+q}),
d(x_{n+1},x_{n+q+1})< \veps+\de/2< \veps+\de, \\
d(x_n,x_{n+1}), d(x_{n+q},x_{n+q+1})< \de/4k<
\de< \veps+\de.
\earr
$$
This, along with the triangular inequality, gives us
$$
d(x_n,x_{n+q+1})\le 
d(x_n,x_{n+q})+d(x_{n+q},x_{n+q+1})< 
\veps+\de/2+\de/4k< \veps+\de;
$$
wherefrom 
$B_1(x_n,x_{n+q})< \veps+\de$; 
so that (by the diameter boundedness property),
$(0<)\ G(x_n,x_{n+q})< \veps+\de$.
Taking the Meir-Keeler
$(d,\calr;G)$-contractive assumption 
imposed upon $T$ into account, gives
$$
d(x_{n+1},x_{n+q+1})
=d(Tx_n,Tx_{n+q})\le \veps;
$$
so that, by the triangular inequality
(and (\ref{404}) again) 
$$
\barr{l}
d(x_n,x_{n+p+1})\le 
d(x_n,x_{n+1})+d(x_{n+1},x_{n+q+1})+
d(x_{n+q+1},x_{n+p+1}) \\
< \veps+\de/4k+k\de/4k\le \veps+\de/4+\de/4
=\veps+\de/2;
\earr
$$
and our claim follows.

{\bf Part 5.}
As  $X$ is $(a-o,d)$-complete,
$x_n\dtends z$, for some 
(uniquely determined) $z\in X$.
If there exists a sequence of ranks $(i(n); n\ge 0)$ 
with [$i(n)\to \infty$ as $n\to \infty$] such that  
$x_{i(n)}=z$ (hence, $x_{i(n)+1}=Tz$) for all $n$, then,
as $(x_{i(n)+1}; n\ge 0)$ is a subsequence of 
$(x_n; n\ge 0)$, one gets $z=Tz$.
So, in the following, we may assume that
the opposite alternative is true:
\ben
\item[] (d04)\ \ 
$\exists h\ge 0$:\ $n\ge h$ $\limpl$ $x_n\ne z$. 
\een
There are several cases to discuss.

{\bf Case 5a.}
Suppose that $T$ 
is $(a-o,d)$-continuous.
Then 
$y_n:=Tx_n\dtends Tz$ as $n\to \oo$.
On the other hand, $(y_n=x_{n+1}; n\ge 0)$ is a subsequence of $(x_n)$;
whence $y_n \dtends z$;
and this yields (as $d$ is sufficient), $z=Tz$.

{\bf Case 5b.}
Suppose that
$\calr$ is $(a-o,d)$-almost-selfclosed.
Put, for simplicity reasons, $b:=d(z,Tz)$.
By definition, 
there exists a subsequence $(u_n:=x_{i(n)}; n\ge 0)$ of 
$(x_n; n\ge 0)$, such that 
$u_n\calr z$, $\forall n$.
Note that, as $\lim_n i(n)= \oo$, one may arrange for
$i(n)\ge n$, $\forall n$; so that, from (d04),
\beq \label{406}
\mbox{
$\forall n\ge h$:\ 
[$i(n)\ge h$;\ whence (see above),\ $u_n\wt \calr z$].
}
\eeq
This, along with 
$(Tu_n=x_{i(n)+1}; n\ge 0)$ being as well
a subsequence of $(x_n; n\ge 0)$, gives 
(via (\ref{402}) and Lemma \ref{le2})
\beq \label{407}
\barr{l}
A_1(u_n,z)=d(u_n,z)\to 0,\  d(Tu_n,z)\to 0,\\
d(u_n,Tu_n)\to 0,\
d(u_n,Tz)\to b,\
d(Tu_n,Tz)\to b;
\earr
\eeq
whence (by definition)  
\beq \label{408}
A_2(u_n,z),A_4(u_n,z)\to b/2,\ 
A_3(u_n,z),B_1(u_n,z)\to b.
\eeq

We now show that the assumption 
$z\ne Tz$ (i.e.: $b> 0$)
yields a contradiction.
Two alternatives must be treated.

{\bf Alter 1.}
Suppose that $G\in \calg_1$.
By the 
Meir-Keeler contractive condition,
$$
d(Tu_n,Tz)< G(u_n,z)\le B_1(u_n,z),\ \forall n\ge h;
$$
so that, combining with the preceding relations,
$G(u_n,z)\to b$.
This, along with (\ref{407})+(\ref{408}),
is impossible for any $G\in \calg_1$;
whence, $z=Tz$.

{\bf Alter 2.}
Suppose that $G\in \calg_2$.
The above convergence properties of $(u_n; n\ge 0)$
tell us that,
for a certain rank $n(b)\ge h$, we must have
$$
d(u_n,Tu_{n}), d(u_n,z), d(Tu_{n},z)< b/2,\ \forall n\ge n(b).
$$
This, by the $d$-Lipschitz property of $d(.,.)$, gives 
$$
|d(u_n,Tz)-b|\le d(u_n,z)< b/2, \forall n\ge n(b),
$$
wherefrom: 
$b/2< d(u_n,Tz)< 3b/2,\ \forall n\ge n(b)$.
Combining these, yields
\beq \label{409}
G(u_n,z)=b,\ \forall n\ge n(b),\  \forall G\in \calg_2.
\eeq
Two sub-cases are now under discussion.

{\bf Alter 2a.}
Suppose that $T$ is 
$(d,\calr;G,\vphi)$-contractive,
for a certain Meir-Keeler
admissible function $\vphi\in \calf(re)(R_+)$.
(The case $G\in \calg_1$ was already clarified
in a preceding step).
By (\ref{409}) and this contractive property,
$$
d(Tu_{n},Tz)\le \vphi(b),\ \forall n\ge n(b).
$$
Passing to limit gives
(by (\ref{407}) above), $b\le \vphi(b)$; 
contradiction; hence, $z=Tz$.

{\bf Alter 2b.}
Suppose that
$T$ is 
$(d,\calr;G,(\psi,\vphi))$-contractive,
for a certain pair $(\psi,\vphi)$ of
weak generalized altering functions in $\calf(R_+)$.
(As before, the case $G\in \calg_1$ 
is clear, by a preceding step).
From this contractive condition,
$$
\psi(d(Tu_n,Tz))\le \psi(G(u_n,z))-\vphi(G(u_n,z)),
\forall n\ge n(b);
$$
or, equivalently (combining with (\ref{409}) above)
$$
0< \vphi(b)\le \psi(b)-\psi(d(Tu_n,Tz)),\
\forall n\ge n(b).
$$
Note that, as a consequence,
$d(Tu_n,Tz)< b$,\ $\forall n\ge n(b)$.
Passing to limit as $n\to \oo$
and taking (\ref{407}) into account, yields 
$\vphi(b)\le \psi(b)-\psi(b-0)$.
This, however, contradicts the choice of
$(\psi,\vphi)$;
so that, $z=Tz$.
The proof is complete.
\eproof

In particular, when $\calr$ is 
transitive, this result is comparable with
the one in 
Turinici \cite{turinici-2011}.
Note that, further extensions 
of these facts are possible,
in the realm of
triangular symmetric spaces,
taken as in 
Hicks and Rhoades \cite{hicks-rhoades-1999};
or, in the setting of
partial metric spaces,
introduced under the lines in 
Matthews \cite{matthews-1994}; 
we shall discuss them elsewhere.

\section{Further aspects}
\setcounter{equation}{0}

Let in the following
$(X,d,\calr)$ be a relational metric space;
and $T$ be a selfmap of $X$.
Technically speaking, 
Theorem \ref{t2} we just exposed 
consists of three sub-statements;
according to the alternatives
of our main result we already listed.
For both practical and theoretical reasons,
it would be useful to evidentiate them; 
further aspects involving the obtained 
facts are also discussed.

Before doing this, 
let us remark that, the condition
\ben 
\item[] (e01)\ \ 
$\calr$ is locally finitely transitive and non-identical
\een
appears as a particular case of (d01).
On the other hand, (d01) is not deductible from
(e01). In fact, 
(d01) has nothing to do with 
the points of 
\ben
\item[] (e02)\ \ 
$X^c(T,\calr):=X\sm X(T,\calr)=\{x\in X; (x,Tx)\notin \calr\}$.
\een
So, even if the restriction of
$\calr$  to 
$X^c(T,\calr)$ is arbitrarily taken, 
(d01) may hold.
On the other hand, (e01) cannot hold 
whenever $X^c(T,\calr)$ 
admits a denumerable subset $Y$ 
such that the 
restriction of $\calr$ to $Y$ is finitely transitive;
and this proves our assertion.

We may now pass to the
particular cases of Theorem \ref{t2}
with practical interest.

{\bf Case 1.}
As a direct consequence of Theorem \ref{t2}, 
we get

\btheorem \label{t3}
Assume that
$T$ is 
$\calr$-semi-progressive,
$\calr$-increasing,
and
Meir-Keeler $(d,\calr;G)$-contractive,
for some $G\in \calg$.
In addition, 
let 
$\calr$ be finitely semi-recurrent non-identical,
$X$ be $(a-o,d)$-complete,
and one of the conditions below holds:

{\bf i1)}
$T$ is $(a-o,d)$-continuous

{\bf i2)}
$\calr$ is $(a-o,d)$-almost-selfclosed
and $G\in \calg_1:=\{A_1,B_2,B_4,C_1\}$.

\n
Then $T$ is a globally strong Picard operator 
(modulo $(d,\calr)$).
\etheorem

The following particular cases 
of this result are to be noted.

{\bf 1-1)}
Let $\sig(.)$ 
be a function in $\calf(X\times X,R_+)$;
and $\cals$
denote the associated relation:
[$x\cals y$ iff $\sig(x,y)\ge 1$].
Then, if we take 
$\calr:=\cals$
and $G=A_1$,
the alternative {\bf i1)} of
Theorem \ref{t3} 
includes 
the related statement in
Berzig and Rus \cite{berzig-rus-2013}.
By the previous remark,
this inclusion is --
at least from a technical viewpoint --
effective;
but, from a logical perspective,
it is possible that the
converse inclusion be also true.
Finally, the alternative {\bf i2)} of 
Theorem \ref{t3}
seems to be new.

{\bf 1-2)}
Suppose that $\calr=X\times X$ 
(i.e.: $\calr$ is the {\it trivial relation} over $X$).
Then, 
Theorem \ref{t3} is comparable with the main
results in
Wlodarczyk and Plebaniak 
\cite{wlodarczyk-plebaniak-2008,
wlodarczyk-plebaniak-2012a, 
wlodarczyk-plebaniak-2012b, 
wlodarczyk-plebaniak-2013},
based on contractive type conditions involving generalized pseudodistances.
However, none of these is reducible to the remaining ones;
we do not give details.

{\bf Case 2.}
As another consequence of Theorem \ref{t2}, we have
the following statement (with practical value):

\btheorem \label{t4}
Assume that $T$ is 
$\calr$-semi-progressive,
$\calr$-increasing,
and
$(d,\calr;G,\vphi)$-contractive,
for some $G\in \calg$
and a certain Meir-Keeler admissible function $\vphi\in \calf(re)(R_+)$.
In addition, 
let 
$\calr$ be finitely semi-recurrent non-identical,
$X$ be $(a-o,d)$-complete,
and one of the conditions below holds:

{\bf j1)}
$T$ is $(a-o,d)$-continuous

{\bf j2)}
$\calr$ is $(a-o,d)$-almost-selfclosed.

\n
Then $T$ is a globally strong Picard operator 
(modulo $(d,\calr)$).
\etheorem

The following particular cases 
of this result are to be noted.

{\bf 2-1)}
Suppose that 
$\calr=X\times X$ (=the trivial relation over $X$)
and $G=A_1$.
Then, Theorem \ref{t4} is comparable with the
main results in
Wlodarczyk et al 
\cite{wlodarczyk-plebaniak-obczynski-2007a,
wlodarczyk-plebaniak-obczynski-2007b},
based on contractive type conditions like
\ben
\item[] (e03)\ \ 
$\diam(T(Y))\le \vphi(\diam(Y))$, for all $Y\in CB(X)$.
\een
[Here, $CB(X)$ is the class of all (nonempty) closed bounded 
subsets of $X$].
Clearly, this condition is stronger than the 
one we already used in Theorem \ref{t4}.
On the other hand, (e03) is written in terms of
generalized pseudodistances.
Hence, direct inclusions between
these results are not in general available;
we do not give details.

{\bf 2-2)}
Suppose that 
$\calr=X\times X$;
and $\vphi\in \calf(re)(R_+)$ is
BWM-admissible
[i.e.: it is either 
Boyd-Wong admissible or Matkowski admissible].
Then, if $G=A_1$, 
Theorem \ref{t4} includes the
Boyd-Wong's result \cite{boyd-wong-1969}
when $\vphi$ is Boyd-Wong admissible;
and, respectively, 
the Matkowski's result \cite{matkowski-1975}
when $\vphi$ is Matkowski admissible.
Moreover, when $G=C_2$,
Theorem \ref{t4} includes
the result in
Leader \cite{leader-1979}.

{\bf 2-3)}
Suppose that $\calr$ is an order on $X$.
Then, 
Theorem \ref{t4} includes the 
results in 
Agarwal et al \cite{agarwal-el-gebeily-o-regan-2008};
see also 
O'Regan and Petru\c{s}el \cite{o-regan-petrusel-2008}.

{\bf Case 3.}
As a final consequence of Theorem \ref{t2}, we have

\btheorem \label{t5}
Assume that the selfmap $T$ is 
$\calr$-semi-progressive,
$\calr$-increasing,
and
$(d,\calr;G,(\psi,\vphi))$-contractive,
for a certain $G\in \calg$
and some pair $(\psi,\vphi)$ of
generalized altering functions in $\calf(R_+)$.
In addition, 
let 
$\calr$ be finitely semi-recurrent non-identical,
$X$ be $(a-o,d)$-complete,
and one of the conditions below holds:

{\bf k1)}
$T$ is $(a-o,d)$-continuous

{\bf k2)}
$\calr$ is $(a-o,d)$-almost-selfclosed.

\n
Then $T$ is a globally strong Picard operator 
(modulo $(d,\calr)$).
\etheorem

The following particular cases 
of this result are to be noted.

{\bf 3-1)}
Let $\al(.),\be(.)$ 
be a couple of functions in $\calf(X\times X,R_+)$;
and $\cala,\calb$
stand for the associated relations
\ben
\item[]
$x\cala y$ iff $\al(x,y)\le 1$;\
$x\calb y$ iff $\be(x,y)\ge 1$.
\een
Then, if we take 
$\calr:=\cala\cap \calb$
and $G\in \calg$,
this result includes (cf. Lemma \ref{le1}) the one in
Berzig et al \cite{berzig-karapinar-roldan-2013},
based on global contractive conditions like
\ben
\item[] 
$\psi(d(Tx,Ty))\le \al(x,y)\psi(d(x,y))-\beta(x,y)\vphi(d(x,y))$,\
$\forall x,y\in X$;
\een
referred to as: 
$T$ is {\it $(\al\psi,\be\vphi)$-contractive}.
In particular, when 
$G=A_1$,
this last result reduces to the one in
Karapinar and Berzig \cite{karapinar-berzig-2013};
which, in turn, extends the one due to 
Samet et al \cite{samet-vetro-vetro-2012};
hence, so does Theorem \ref{t5} above.

{\bf 3-2)}
Let $(Y,d)$ be a metric space;
and $T$ be a selfmap of $Y$.
Given $p\ge 2$, let $\{A_1,...,A_p\}$
be a finite system of closed subsets of
$Y$ with
\ben
\item[] (e04)\ \ 
$T(A_i)\incl A_{i+1}$,\ 
for all $i\in \{1,...,p\}$\ (where $A_{p+1}=A_1$).
\een
Define a relation $\calr$ over $Y$ as
\ben
\item[] (e05)\ \ 
$\calr=(A_1\times A_2)\cup ... \cup(A_p\times A_{p+1})$;
\een
then, put $X=A_1\cup ... \cup A_p$.
Clearly, $T$ is a selfmap of $X$; and the relation
$\calr$ is $p$-semi-recurrent 
at each orbital admissible point of $X(T,\calr)$.
The corresponding version of Theorem \ref{t5}
includes the related statement in
Berzig et al \cite{berzig-karapinar-roldan-2013}.

It is to be stressed that, 
this last construction may be also attached to 
the setting of {\bf Case 2}.
Then, the corresponding version of Theorem \ref{t4}
extends in a direct way some basic results in
Kirk et al \cite{kirk-srinivasan-veeramani-2003}.

Finally, we should remark  that, 
none of these particular theorems 
may be viewed as a genuine extension for 
the fixed point statement due to
Samet and Turinici \cite{samet-turinici-2012};
because, in the quoted paper,
$\calr$ is not subjected to
any kind of 
(local or global)
transitive type requirements.
Further aspects 
(involving the same general setting) 
may be found in 
Berzig \cite{berzig-2013}.


\end{document}